\documentclass[aps,nofootinbib]{revtex4}

\usepackage{mathpazo}

\usepackage{amsmath}
\usepackage{amsfonts,amssymb,amsthm, bbm,braket, color}
\usepackage{graphicx}   
\usepackage{subfigure}  
\usepackage{amsbsy} 
\usepackage[bold]{hhtensor} 
\usepackage{natbib}
\usepackage{algpseudocode}

\usepackage{multirow}

\newcommand\Tr{\mathrm{Tr}}
\DeclareMathOperator*{\argmin}{arg\,min}
\DeclareMathOperator*{\argmax}{arg\,max}
\def\Id{1\!\mathrm{l}}

\usepackage{hyperref} 

\hypersetup{
  colorlinks   = true, 
  urlcolor     = blue, 
  linkcolor    = blue, 
  citecolor   =  red 
}

\usepackage{algorithm} 

\begin{document}

\title{Bayes estimator for multinomial parameters and Bhattacharyya distances}

\author{Christopher Ferrie}
\affiliation{Centre for Engineered Quantum Systems, School of Physics, The University of Sydney, Sydney, NSW, Australia}

\author{Robin Blume-Kohout}
\affiliation{Sandia National Laboratories, Albuquerque, New Mexico, 87185}

\date{\today}


\begin{abstract}
We derive the Bayes estimator for the parameters of a multinomial distribution under two loss functions ($1-B$ and $1-B^2$) that are based on the Bhattacharyya coefficient $B(\vec{p},\vec{q}) = \sum{\sqrt{p_kq_k}}$.  We formulate a non-commutative generalization relevant to quantum probability theory as an open problem.  As an example application, we use our solution to find minimax estimators for a binomial parameter under Bhattacharyya loss ($1-B^2$).
\end{abstract}


\maketitle


\section{Introduction}

In statistical decision theory, a Bayes estimator for a prior $\mu$ is a point estimator whose expected risk (with respect to $\mu$) is as small as possible \cite{Berger}.  Bayes estimators have utility beyond strictly Bayesian statistics; they also provide lower bounds on \emph{minimax} risk (a frequentist concept), and the greatest lower bound (That is, the maximum Bayes risk over all priors) coincides with the minimax risk.  This duality enables efficient numerical algorithms to find minimax estimators by optimizing over Bayesian priors, rather than over the computationally intractable space of estimators \cite{Kempthorne1987Numerical}.

Different loss functions produce different risk functions, and the Bayes estimator is sensitive to this choice.  For a certain class of risks known as Bregman divergences the Bayes estimator is always the mean of the posterior distribution \cite{banerjee}.  For the multinomial model that we consider, both mean-squared-error and Kullback-Leibler divergence (relative entropy) are Bregman divergences.  But not every useful risk function has this property.  We consider two loss functions defined through the Bhattacharyya coefficient $B(\vec{p},\vec{q})$, a Fisher-adjusted measure of distinguishability between probability distributions that finds use in machine learning (as a metric on probability distributions \cite{mlBhat}) and in quantum information theory (as \emph{fidelity} \cite{Wootters, Jozsa}, the most commonly used measure of quantum state distinguishability in both theoretical and experimental studies).

In this paper, we find the Bayes estimators for $1-B$ and $1-B^2$ by mapping their respective optimizations to simple linear programs.  We also state, but do not solve, the non-commutative generalization of the problem, which has applications in quantum information theory. 
Finally, we demonstrate the utility of the Bayes estimator by finding the minimax estimator for a binomial parameter.

\section{Problem set-up, review and choice of loss functions \label{sec:1}}

Suppose we have a coin, a $K$-sided die, or any process that generates i.i.d. samples from a finite set of size $K$.  We wish to estimate the probabilities $\{p_0,\ldots p_{K-1}\}$.  They form a parameter vector $\vec{p}$ which belongs to the $K$-simplex
\begin{equation}
\triangle_K = \left\{(p_0, p_1\ldots,p_{K-1})\bigg| p_k\geq 0 \;\forall k, \sum_k p_k =1\right\}.
\end{equation}
Estimation requires some data.  From the data, we wish to produce an estimate of $\vec{p}$, call it $\hat{\vec{p}}$.  Although the data are \emph{most} commonly obtained by repeated sampling directly from $\vec{p}$, other experiments are possible.  We might have ``noisy'' data (drawn from an ancillary distribution $\vec{q} = \Sigma\vec{p}$, where the stochastic matrix $\Sigma$ is known \cite{NoisyCoinPaper}), or side information.  Happily, the Bayes analysis turns out to be independent of how the data were obtained.  Instead, it depends only on the posterior probability density $\Pr(\vec{p}|\texttt{data})$---which is to say that, given the posterior, the expected risk is independent of the likelihood $\Pr(\texttt{data}|\vec{p})$.

The quality of an estimate $\hat{\vec{p}}$ is formalized by a loss function $L(\vec{p},\hat{\vec{p}})$ that quantifies how bad it is to report an estimate $\hat{\vec{p}}$ when the truth is $\vec{p}$.  The best estimate should be $\vec{p}$ itself, so $L(\vec{p},\vec{p})\leq L(\vec{p},\hat{\vec{p}})$ for all $\hat{\vec{p}}$.  A ``good'' estimator is one that minimizes expected loss, or \emph{risk}:
\begin{equation}
r_{\vec{p}} = \mathbb E_{\texttt{data}}[L(\vec{p},\hat{\vec{p}}(\texttt{data})].
\end{equation}
The risk is a function of the underlying $\vec{p}$, and it is generally impossible to minimize risk for \emph{all} $\vec{p}$ simultaneously.  The Bayesian solution is to consider (and minimize) the average risk with respect to a given prior $\mu(\vec{p})$.  This is called the \emph{Bayes risk} of the estimator $\hat{\vec{p}}$ with respect to $\mu$:
\begin{equation}\label{bayes risk}
r(\hat{\vec{p}}) = \mathbb E_{\vec{p},\texttt{data}}[L(\vec{p},\hat{\vec{p}}(\texttt{data})].
\end{equation}
An estimator $\hat{\vec{p}}_{\rm B}$ that minimizes the Bayes risk is called a \emph{Bayes estimator} for $\mu(\vec{p})$.  The Bayes estimator can also be defined (under normal circumstances; \cite{Lehmann1998Theory}) explicitly for each dataset as the minimizer of \emph{posterior} risk.  This identity is extremely useful, because it means that a Bayes estimator can be constructed explicitly for any given value of \texttt{data} without considering others.

Different loss functions yield different Bayes estimators.  For example, for quadratic loss functions of the form
\begin{equation}
L(\vec{p},\hat{\vec{p}}) = (\vec{p} - \hat{\vec{p}})^\text{T} A (\vec{p} - \hat{\vec{p}}),
\end{equation}
where $A> 0$ is a positive definite matrix, the Bayes estimator is the mean of the posterior distribution \cite{Lehmann1998Theory}
\begin{equation}
\hat{\vec{p}}_{\rm B}(\texttt{data})= \mathbb E_{\vec{p}|\texttt{data}}[\vec{p}].
\end{equation}
Another important loss function is the \emph{relative entropy} or Kullback-Leibler divergence:
\begin{equation}
L(\vec{p},\hat{\vec{p}}) = \sum_k p_k \log \frac{p_k}{\hat p_k}.
\end{equation}
Here, again, the Bayes estimator is the posterior mean \cite{Aitchison}.  This is not an accident; it holds true whenever the loss is a \emph{Bregman divergence} \cite{banerjee}.  Quadratic loss and Kullback-Leibler loss are both Bregman divergences.

Unfortunately, not every useful and important loss function is a Bregman divergence.  The loss functions we consider here are not Bregman divergences, and their Bayes estimators are not the posterior mean.  They are based on the Bhattacharyya coefficient \cite{bhat},
\begin{equation}
B(\vec{p},\hat{\vec{p}}) = \sum_k \sqrt{p_k \hat p_k},
\end{equation}
which quantifies the statistical indistinguishability or \emph{fidelity} of two distributions.
In different contexts, both $L_1 = 1-B$ and $L_2 = 1-B^2$ are useful loss functions, and so we find Bayes estimators for both of them.  In the limit of small deviations ($\hat{\vec{p}}$ is close to $\vec{p}$), they coincide up to a factor of 2.  Both are Fisher-adjusted, meaning that their 2nd order series expansions are proportional to the Fisher metric.  However, $L_1$ and $L_2$ have different global behavior, and therefore different Bayes estimators.

Since the minimizer of $1-B^p$ is the maximizer of $B^p$, we define two Bayes estimators as follows:
\begin{align}	
	\hat {\vec p}_{1} &= \underset{\hat {\vec p}}{\rm argmax} \;\mathbb E_{\vec p}[B(\vec p, \hat{\vec p})], \label{p_1}\\
	\hat {\vec p}_{2} &= \underset{\hat {\vec p}}{\rm argmax} \;\mathbb E_{\vec p}[B^2(\vec p, \hat{\vec p})],\label{p_2}
\end{align}
subject to $\hat{\vec p} \in \triangle_K$. Note that, since the Bayes estimator depends only on the posterior, with no explicit dependence on the prior or the data, we hereafter drop the conditional notation on the data. That is, $\mathbb E_{\vec p}$ means expectation with respect to an arbitrary distribution of $\vec{p}$.

\section{Bayes estimators \label{sec:3}}

Let us now derive expressions for the Bayes estimators defined in Equations \ref{p_1} and \ref{p_2}.  A very useful way of writing $B$ is to define $\sqrt{\vec p}$ as the \emph{element-wise} square root of $\vec p$, then write:
\begin{equation}
B(\vec{p},\hat{\vec{p}}) = \sqrt{\vec p}\cdot \sqrt{{\hat{\vec p}}}.
\end{equation}
The posterior risk involves an average of $B$ over $\vec{p}$.  Both this average and the dot product in $B$ are \emph{linear} in $\sqrt{\vec p}$, which means that they commute, so:
\begin{equation}
\mathbb E_{\vec p}[B(\vec{p},\hat{\vec{p}})] = \mathbb E_{\vec p}[\sqrt{\vec p}]\cdot \sqrt{{\hat{\vec p}}}.
\end{equation}
This also implies the normalization
\begin{equation}
	1 = \sum_k p_k = \sum_k \sqrt{p_k}\sqrt{p_k} = \|\sqrt{\vec p}\|^2.
\end{equation}
Thus, the problem in \autoref{p_1} becomes
\begin{equation}\label{problem 1}
\begin{aligned}
& {\text{maximize}}
& & \mathbb E_{\vec p}[\sqrt{\vec p}]\cdot \sqrt{{\hat{\vec p}}} \\
& \text{subject to}
& & \sqrt{\hat{\vec p}}\cdot \sqrt{\hat{\vec p}} = 1, \\
& & & \sqrt{\hat{\vec p}} \geq 0.
\end{aligned}
\end{equation}
This is a textbook \cite{boyd} convex optimization problem, and the solution is
\begin{equation}
	\sqrt{\hat {\vec p}_1} = \frac{\mathbb E_{\vec p}[\sqrt{\vec p}]}
	                              {\| \mathbb E_{\vec p}[\sqrt{\vec p}] \|},
\end{equation}
implying
\begin{equation}
	\hat {\vec p}_1 = \frac{\mathbb E_{\vec p}[\sqrt{\vec p}]^2}
	                       { (\mathbb E_{\vec p}[\sqrt{\vec p}])\cdot(\mathbb E_{\vec p}[\sqrt{\vec p}])},
\label{eq:p1sol}
\end{equation}
where the square in the numerator is \emph{element-wise}.

Minimizing $1-B^2$, as in \autoref{p_2}, requires a bit more work, but we can again reduce it to a linear program.  We start by defining
\begin{equation}
	P = \sqrt{\vec p}\sqrt{\vec p}^{\rm T},
\end{equation}
the projection operator onto $\sqrt{\vec p}$, and similarly for $\hat P$. The normalization condition on $P$ is
\begin{equation}
	\Tr[P] = 1.
\end{equation}
In terms of these objects, the squared Bhattacharyya coefficient is
\begin{equation}
B^2(\vec{p},\hat{\vec{p}}) = \Tr[P\hat P].
\end{equation}
This expression is now linear in $P$, which allows us to apply the trick of commuting the expectation through:
\begin{equation}
\mathbb E_{\vec p}[B^2(\vec{p},\hat{\vec{p}})] = \Tr[\mathbb E_{\vec p}[P]\hat P].
\end{equation}
The matrix elements of the expectation on the RHS can be written out explicitly as
\begin{equation}
	\left(\mathbb E_{\vec p}[P]\right)_{ij} = \mathbb E_{\vec p}[\sqrt{p_i p_j}].
\end{equation}
Putting everything together, \autoref{p_2} becomes
\begin{equation}\label{problem 2}
\begin{aligned}
& {\text{maximize}}
& & Tr[\mathbb E_{\vec p}[P]\hat P] \\
& \text{subject to}
& & \Tr[\hat P] = 1, \\
& & & \hat P \geq 0,
\end{aligned}
\end{equation}
where $\hat P \geq 0$ means positive semi-definiteness as a matrix.  This is another textbook \cite{boyd} linear semidefinite program.  If we define $\vec a$ as the normalized eigenvector of $\mathbb E_{\vec p}[P]$ with maximal eigenvalue, then the solution to \autoref{problem 2} is
\begin{equation}
	\hat P_2 = \vec a \vec a^{\rm T}.
\end{equation}
Explicitly, as a vector in $\triangle_K$, we have
\begin{equation}
	\hat {\vec p}_2 = \vec a^2,
\end{equation}
where $\vec a = (a_0, a_1, \ldots, a_{K-1})$ is the eigenvector of $\mathbb E_{\vec p}[\sqrt{p_i p_j}]$ with maximal eigenvalue, and the square is element-wise as in Eq. \ref{eq:p1sol}.

\section{Non-commutative generalization: a quantum (open) problem \label{sec:5}}

Quantum mechanics can be described as a non-commutative generalization of probability theory \cite{nc}.  Quantum random variables (That is, the ``state'' of a quantum system) are represented not by $K$-element vectors of probabilities $\vec{p}$, but by $K\times K$ self-adjoint complex matrices or \emph{density matrices}.  In this framework, observable events are also represented by $K\times K$ self-adjoint matrices called \emph{effects}.  States can be sampled or observed in more than one way.  The quantum analogue to ``sampling'' is ``measuring'' the system, and the observer must choose \emph{how} to measure it.  This means specifying a mutually exclusive and exhaustive set of possible \emph{events}, each represented by a positive semidefinite $d\times d$ matrix.  The set of event matrices $\{E_k\}$ is called a \emph{positive operator valued measure} (POVM), and satisfies the condition $\sum_k{E_k} = \Id$.  If an quantum system is state $\rho$ is observed, then the probability of observing the event represented by effect $E$ is given by $\mathrm{Pr}(E|\rho) = \Tr(\rho E)$.  As a result, the law of total probability (``exactly one of the possible events occurs'') manifests itself as the constraint $\Tr\rho=1$, and the non-negativity of probabilities corresponds to $\rho\geq0$.  Thus, the set of quantum states (analogous to the simplex $\Delta$ of valid probability vectors) is
\begin{equation}
\mathcal S = \left\{\rho |\rho\geq 0, \Tr(\rho) = 1\right\}.
\end{equation}
The subset of $S$ containing only diagonal matrices (``classical'' states) is in 1:1 correspondence with $K$-element probability vectors, since $\rho \geq 0$ ensures positivity of the diagonal entries, and $\Tr(\rho) = 1$ ensures their normalization.  Thus, if we consider \emph{only} states that are simultaneously diagonalizable, we can identify $\vec{p}\equiv \rho$ and $\hat{\vec{p}}\equiv \hat\rho$, and the Bhattacharyya coefficient has a simple expression in terms of $\rho$ and $\hat\rho$:
\begin{equation}
B(\vec{p},\hat{\vec{p}}) = \Tr\left(\sqrt{\rho\hat\rho}\right).
\end{equation}
Its square---the $B^2$ risk that we considered previously---is called \emph{classical fidelity} in the quantum mechanics literature.  Hence, our solution to \autoref{problem 2} can be applied directly to estimation of quantum states, whenever the support of the averaging measure (e.g. a posterior distribution) consists entirely of \emph{commuting} quantum states:
\begin{equation}\label{the problem}
\begin{aligned}
& {\text{maximize}}
& & \mathbb E_\rho[\Tr(\sqrt{\rho\hat\rho})^2] \\
& \text{subject to}
& & \Tr(\hat\rho) = 1, \\
&&& \hat\rho \geq 0.
\end{aligned}
\end{equation}

However, quantum states are not usually restricted to be diagonal, which means that in the estimation of quantum states (\emph{quantum tomography} \cite{gill}), the posterior will have support on non-commuting density matrices.  This requires new definitions of loss and risk.  Without getting into the details (see Ref. \cite{Jozsa}), the general definition of \emph{quantum fidelity} between two states is:
\begin{equation}\label{Jfid}
F(\rho,\hat\rho) = \left[\Tr\sqrt{\sqrt\rho \hat\rho \sqrt\rho}\right]^2,
\end{equation}
and the corresponding loss function is \emph{infidelity} or $L = 1-F$.  For co-diagonal states $\hat\rho$ and $\rho$, fidelity coincides with $B^2$.  Thus, if the posterior \emph{does} happen to be supported only on mutually commuting states, then our solution to \autoref{problem 2} gives the Bayes estimator.  In general, however, the problem of finding a Bayes estimator for quantum infidelity is:
\begin{equation}\label{the full problem}
\begin{aligned}
& {\text{maximize}}
& & \mathbb E_{\rho}[F(\rho,\hat\rho)] \\
& \text{subject to}
& & \Tr(\hat\rho) = 1, \\
&&& \hat\rho \geq 0.
\end{aligned}
\end{equation}
This problem appears to be quite difficult (some bounds are given in \cite{kueng}), and we do not attempt to solve it in full generality here.

\section{Example: binomial parameter\label{sec:4}}

\subsection{Canonical beta prior\label{sec:4a}}

\begin{figure}[ht]
\includegraphics[width=0.49\columnwidth]{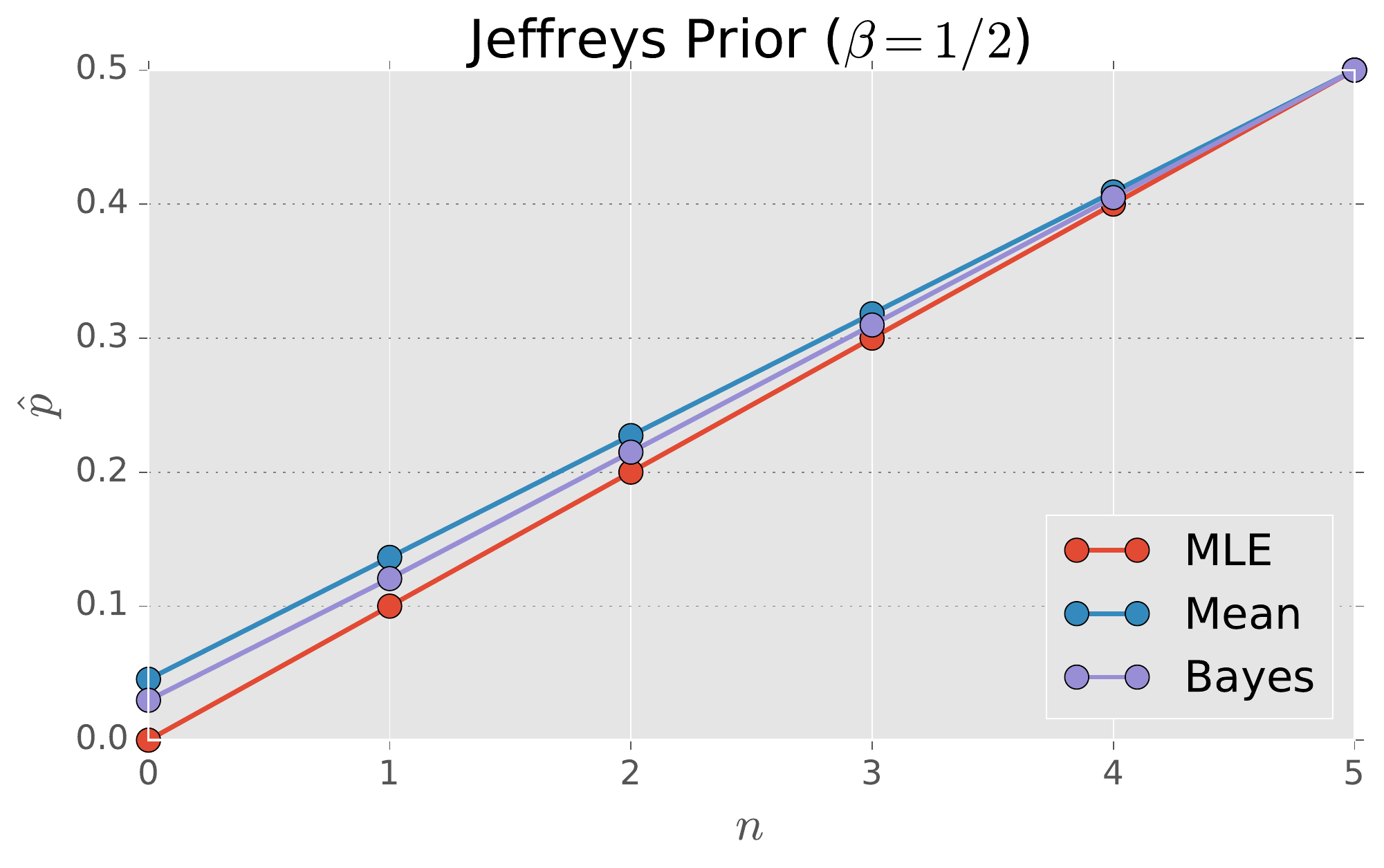}
\includegraphics[width=0.49\columnwidth]{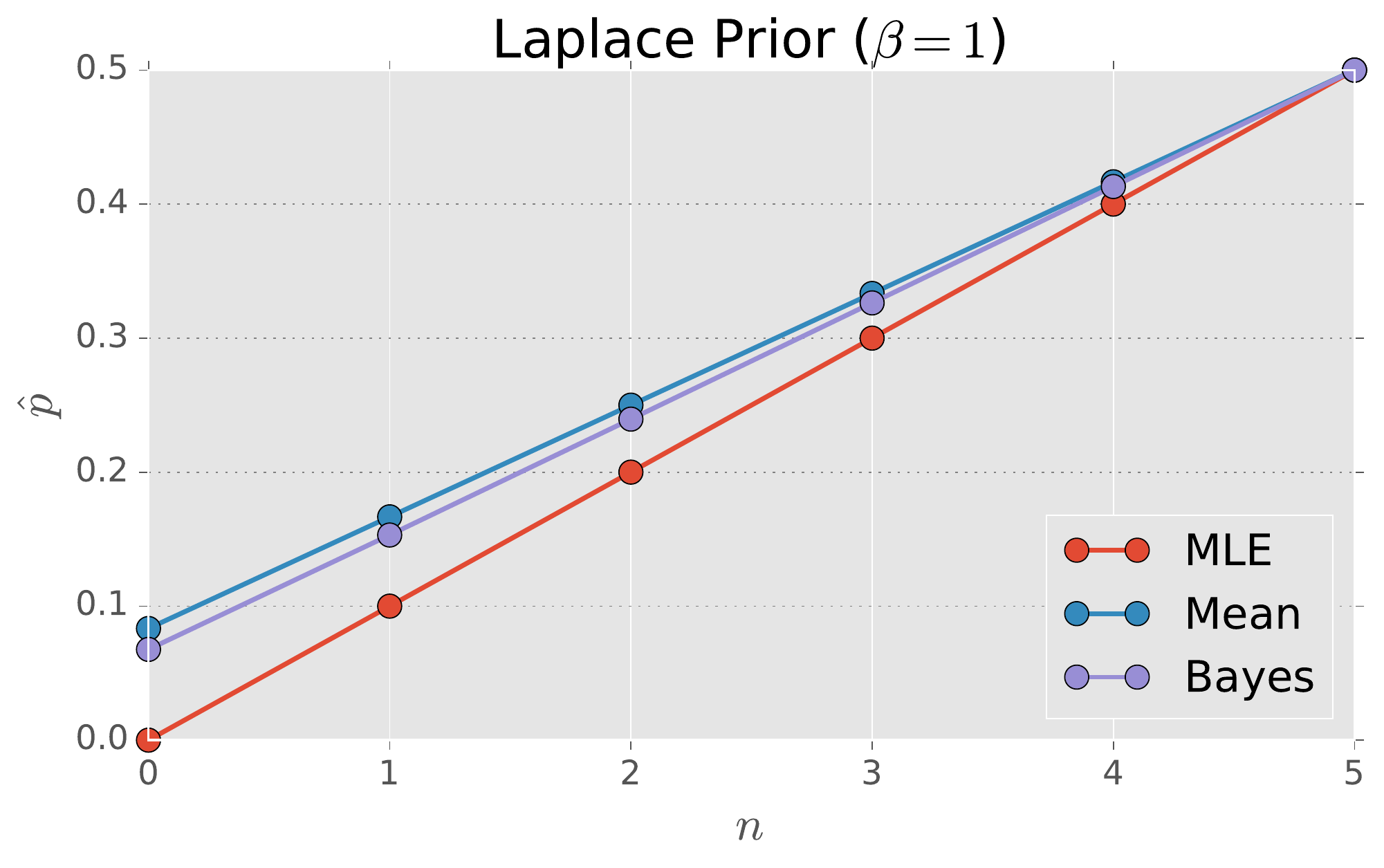}
\caption{\label{fig:estimators_N10} The maximum likelihood estimator (MLE), optimal Bayes estimator and the posterior mean for $N=10$ coin flips.  On the left, the latter two estimators begin with a canonical beta prior with $\beta=1/2$ (Jeffrey's prior).  On the right, the same but with $\beta=1$ (Laplace's uniform prior).  The fidelity performance of these estimators is shown in Figure \ref{fig:fid_N10}.}
\end{figure}

As an interesting and useful application of our main result, we now compute (and examine) the Bayes estimator $\hat{\vec{p}}_2$ for $L_2 = 1-B^2$ in the simple case of a binomial distribution---that is, a directly observed coin toss or Bernoulli process.  Here, $\vec p = (p_0,p_1)$, with $p_1 = 1 - p_0$.  The canonical (conjugate) prior is a beta distribution:
\begin{equation}
\Pr(\vec p) = \frac{p_0^{\beta-1}p_1^{\beta-1}}{{\rm Beta}(\beta,\beta)}d\!\vec{p},
\label{eq:beta}
\end{equation}
where ${\rm Beta}()$ is the beta function.  If the coin is flipped $N$ times, yielding $n$ heads and $N-n$ tails, then the posterior is
\begin{equation}
\Pr(\vec p|n) = \frac{p_0^{n+\beta-1}p_1^{N-n+\beta-1}}{{\rm Beta}(n+\beta,N-n+\beta)}d\!\vec{p}.
\end{equation}
The matrix $ \mathbb E_{\vec{p}}[\sqrt{p_jp_k}]$ is given explicitly in terms of beta functions by
\begin{equation}\label{T mtx}
\left( \begin{array}{cc}
{\rm Beta}(\beta+n+1,\beta+N-n) & {\rm Beta}(\beta+n+1/2,\beta+N-n+1/2) \\
{\rm Beta}(\beta+n+1/2,\beta+N-n+1/2) & {\rm Beta}(\beta+n,\beta+N-n+1)  \end{array} \right)/{\rm Beta}(\beta+n,\beta+N-n).
\end{equation}
To compute the Bayes estimator $\hat{\vec{p}}_2$, we diagonalize this matrix, extract the eigenvector corresponding to the larger eigenvalue, and square its entries.  This can be written down in closed form, but the result is lengthy, messy, and opaque.  
Instead, we illustrate its properties for $N=10$ flips and for conjugate priors with $\beta=1$ (Laplace's prior) and $\beta=1/2$ (Jeffrey's prior), in Fig. \ref{fig:estimators_N10}.  For comparison, we also show (1) the posterior mean estimator, and (2) the maximum likelihood estimator.

The Bayes estimator is generally similar to the posterior mean estimator.  Both ``hedge'' away from $p=0$ and $p=1$, in contrast to the MLE, but the Bayes estimator is more aggressive (less hedged) than the mean.  Since the mean estimator is Bayes for Bregman divergences, this implies that the $1-B^2$ loss function is more forgiving of extreme estimates (That is, ones close to $p=0$ or $p=1$) than any Bregman divergence.

\begin{figure}[ht]
\includegraphics[width=.55\columnwidth]{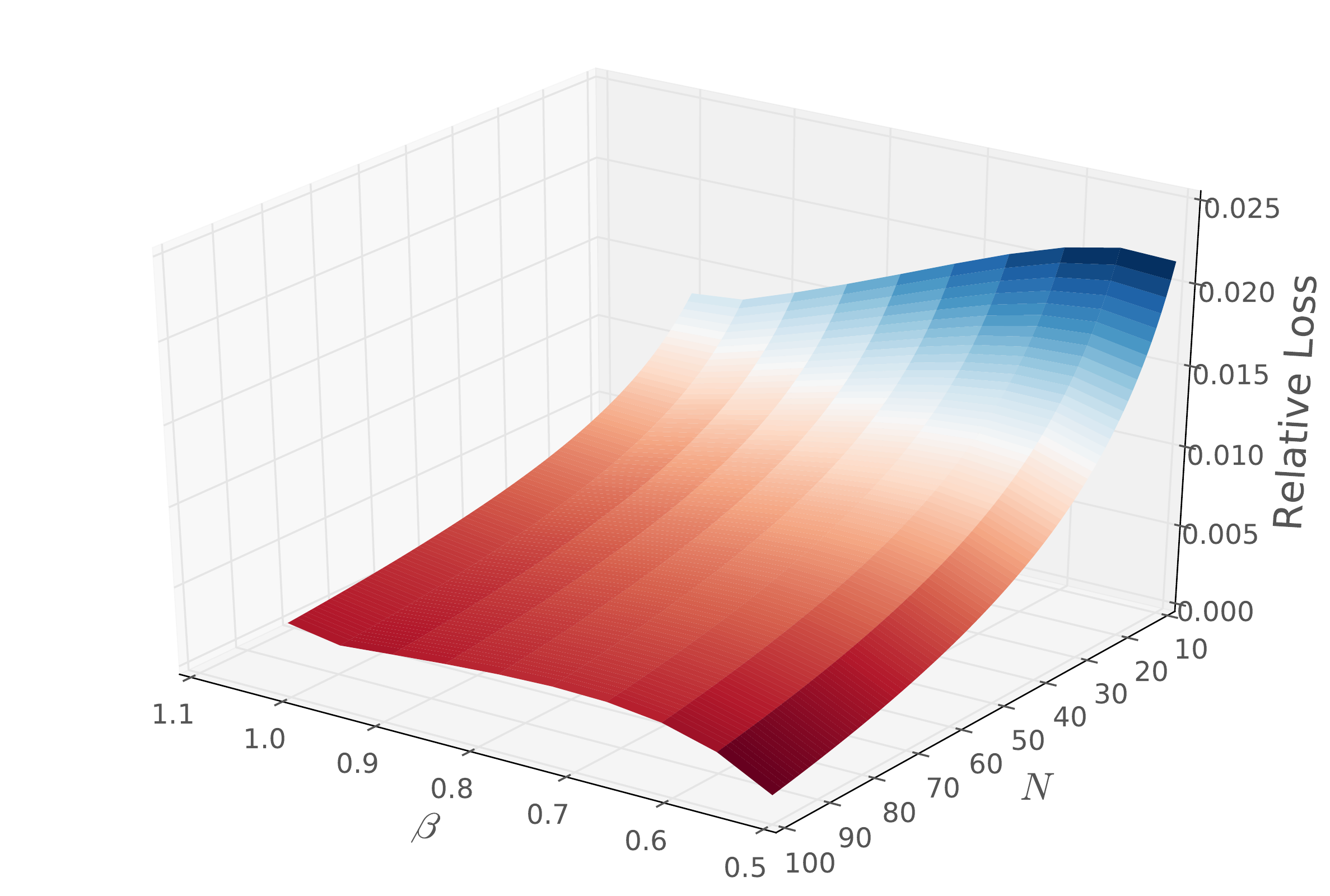}
\caption{\label{fig:reldiff_N10}
The relative difference, $\mathcal R$ in Eq. \eqref{eq:rel}, between the average loss of the Bayes estimator $r_{\rm Bayes}$ and that of the mean estimator $r_{\rm mean}$.}
\label{figreldiff_N10}
\end{figure}

However, although the difference between the two estimators (Bayes and posterior mean) appears significant, the difference in their \emph{performance} is not.  While the Bayes estimator does have lower Bayes risk, the difference is small even for $N=10$.  Figure \ref{figreldiff_N10} shows the \emph{relative suboptimality} of the mean estimator, defined as:
\begin{equation}\label{eq:rel}
\mathcal R = \frac{r_{\rm Bayes}-r_{\rm mean}}{r_{\rm Bayes}},
\end{equation}
where $r_{\rm mean} = \mathbb E_{\vec{p}|n}[1-B^2(\vec{p},\hat{\vec{p}}_{\rm mean}(n))]$ and likewise for $r_{\rm Bayes}$.  The relative suboptimality of the mean estimator---even at $N=10$---is always less than $10^{-3}$, decreases rapidly as a function of $N$ and is largely independent of $\beta$.  

Given that the mean estimator is far easier to calculate (the Bayes estimator requires computing the matrix $ \mathbb E_{\vec{p}}[\sqrt{p_jp_k}]$ and its eigenvectors), we suggest that in practice for most purposes, the mean estimator is a completely satisfactory heuristic.  For small $N$, there is a noticeable difference in the estimate itself (see Fig. \ref{fig:fid_N10}), but \emph{not} in performance (see Fig. \ref{figreldiff_N10}).

\subsection{Least favorable prior\label{sec:4b}}

Frequentist analyses of estimators average the risk \emph{only} over the data (That is with respect to $\Pr(\texttt{data}|\vec{p})$), rather than over the joint distribution $\Pr(\texttt{data},\vec{p})$.  This more restricted average defines \emph{pointwise} risk,
\begin{equation}
R(\vec{p},\hat{\vec{p}}) = \mathbb E_{\texttt{data}|\vec{p}}[L(\vec{p},\hat{\vec{p}})],
\end{equation}
which retains a dependence on $\vec{p}$ that has to be gotten rid of somehow in order to define an ``optimal'' estimator.
Instead of averaging over a prior (to get the Bayes risk in Eq.~\eqref{bayes risk}), the most common frequentist way to remove this dependence is to consider the ``worst case'', and maximize the risk over $\vec{p}$:
\begin{equation}
R_{\text{max}}(\hat{\vec{p}}) = \max_{\vec{p}} R(\vec{p},\hat{\vec{p}}).
\end{equation}

\begin{figure}[ht]
\includegraphics[width=.49\columnwidth]{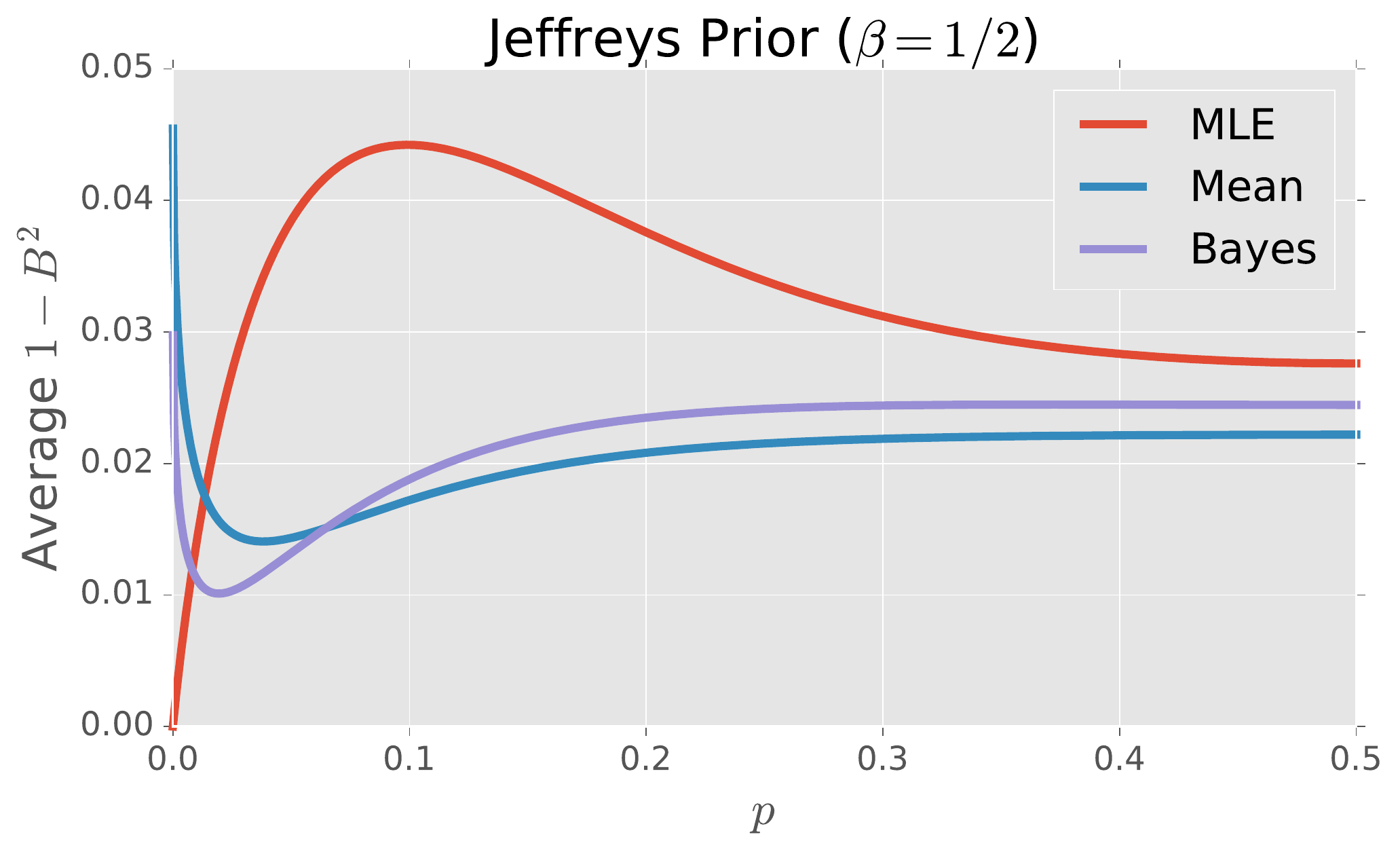}
\includegraphics[width=.49\columnwidth]{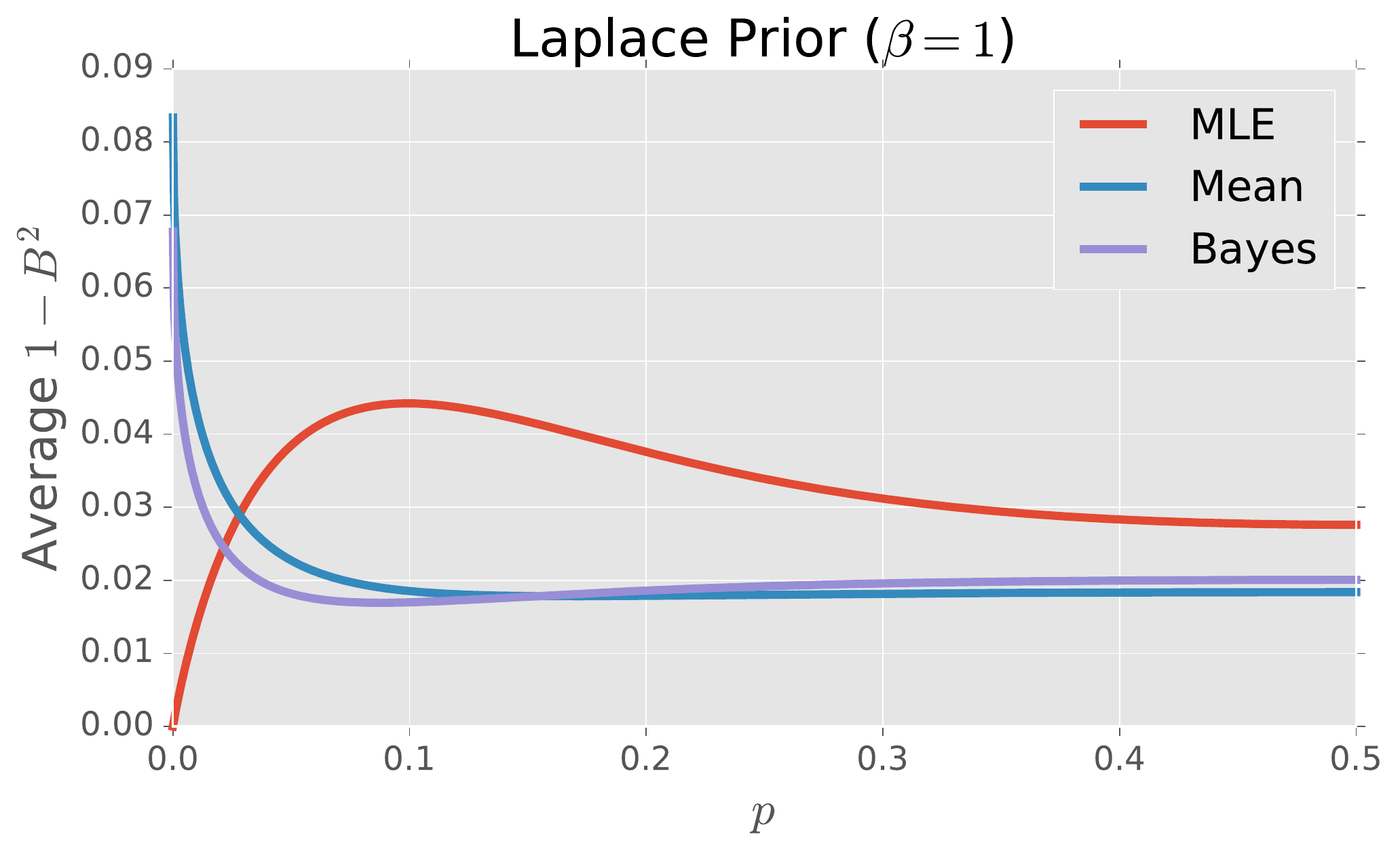}
\caption{\label{fig:fid_N10}
The pointwise risk ($1-B^2$) of the estimators shown in Figure \ref{fig:estimators_N10}, as a function of $p_0$.}
\end{figure}

An estimator that minimizes the maximum risk is called a \emph{minimax} estimator:
\begin{equation}
\hat{\vec{p}}_{\text{minimax}} = \argmin_{\hat{\vec{p}}} R_{\text{max}}(\hat{\vec{p}}).
\end{equation}
The minimax criterion seeks an estimator whose performance is ``pretty good'' for any value of the parameters, without reference to prior probability.  In practice, minimax estimators achieve (approximately) equal risk for all $\vec{p}$.  The Bayes estimators for the binomial parameter that we derived in the previous section are not minimax; Figure \ref{fig:fid_N10} shows that their pointwise risk varies with $p_0$.  However, they're pretty close---the variations are small (e.g., compared with those of the MLE).

However, while the minimax estimator is unique, each prior has a different Bayes estimator.  We can rewrite Eq.~\eqref{bayes risk}, making this explicit, letting $\mu$ denote the prior density:
\begin{equation}
r(\mu,\hat{\vec{p}}) = \mathbb E_{\vec{p}}[R(\vec{p},\hat{\vec{p}})]
\end{equation}
Somewhat remarkably, there \emph{is} (at least under weak regularity conditions \cite{Berger}) a prior whose Bayes estimator is minimax.  The \emph{least favorable prior} is the prior whose Bayes estimator has the highest Bayes risk (out of all priors), and its Bayes estimator is minimax.  Explicitly, the estimator is
\begin{equation}\label{lfp}
\hat{\vec{p}}_{\text{LFP}} = \argmax_{\mu} \min_{\hat{\vec{p}}}r(\mu,\hat{\vec{p}})
\end{equation}
This remarkable relation between Bayesian and frequentist optimality is called \emph{Bayes-minimax duality}:
\begin{equation}
\hat{\vec{p}}_{\text{LFP}} = \hat{\vec{p}}_{\text{minimax}}.
\end{equation}

One useful consequence of Bayes-minimax duality is that it provides a method for finding minimax estimators by searching for the least favorable prior.  While the space of priors is large, it is enormously smaller than the space of \emph{estimators}.  However, this approach is only practical if we have a closed-form expression for Bayes estimators, since this allows fast and efficient computation of the Bayes risk for each prior considered.  Our result enables this sort of analysis for $1-B$ and $1-B^2$ loss functions.

We used our result to construct minimax estimators for the binomial parameter.  Unlike the Bayes estimator, the minimax estimators do not have explicit closed forms; the least favorable priors are not elegant, and require numerical optimization.  We performed two numerical optimizations.  First, we did a \emph{restricted} optimization over conjugate priors (beta distributions, of the form given in Eq. \ref{eq:beta}), to find the minimax value of $\beta$.  Then, we performed an unrestricted optimization over \emph{all} priors to find the true minimax estimator 

In the first (restricted) optimization over conjugate priors, we found that the optimal $\beta$ varied only weakly with $N$, and was given by $\beta\approx 0.44$ for all $N$.  This is extremely consistent with other answers to the general question ``What value of $\beta$ works best?''; a variety of other analyses have yielded answers close to the $\beta=1/2$ that defines Jeffrey's estimator.  Just for reference, we considered what prior would yield the lowest maximum risk if we were to use the posterior mean estimator rather than the Bayes estimator derived in this paper.  For this ad-hoc procedure, the optimal $\beta$ is approximately $\beta\approx0.26$.  This estimator is not minimax in any sense, but it does support our general observation that the Bayes estimator can be approximated by the posterior mean without much damage ($0.26$ and $0.44$ are not very different, and the achieved maximum risks are also not very different).

\begin{figure}[ht]
\includegraphics[width=.49\columnwidth]{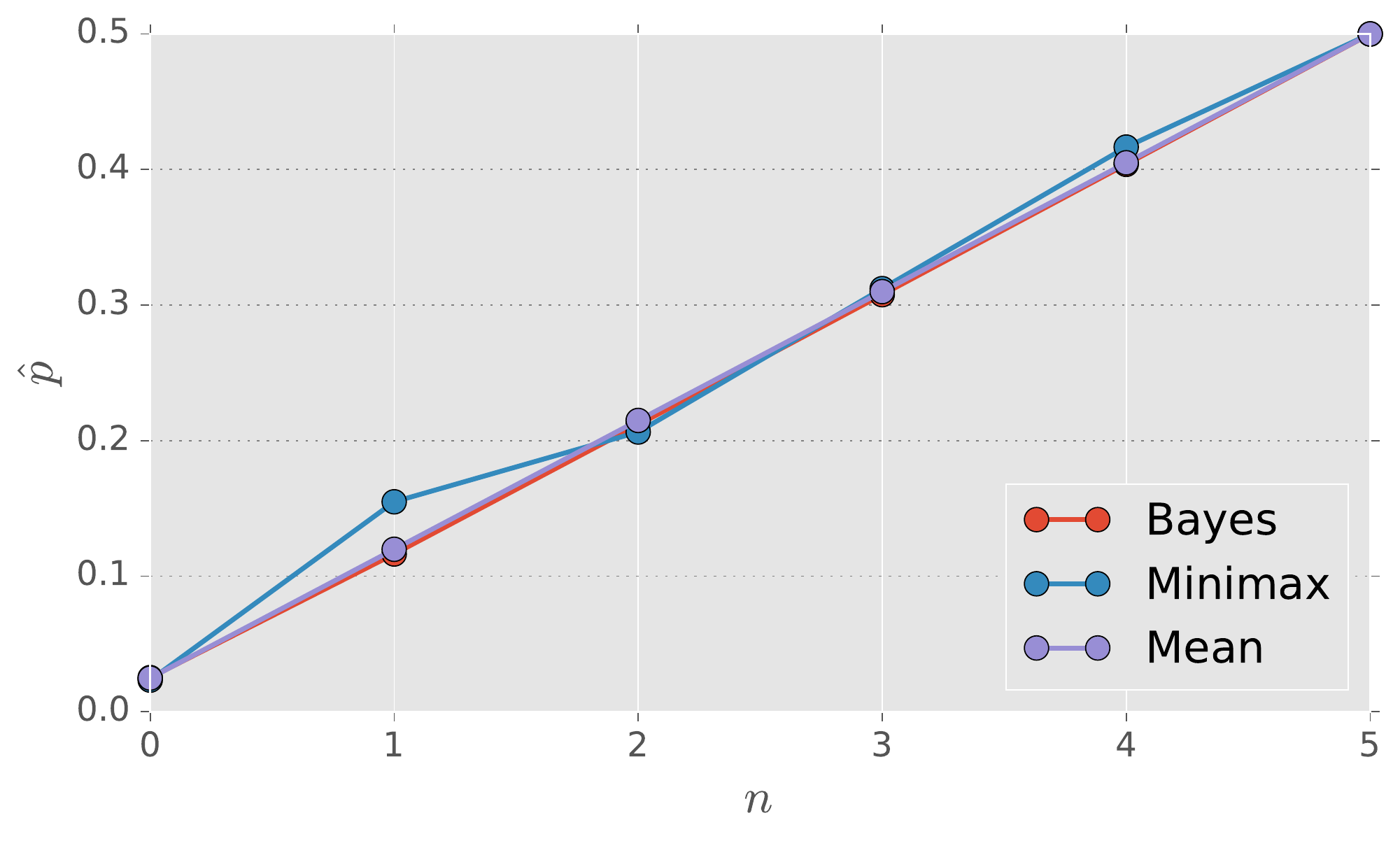}
\includegraphics[width=.49\columnwidth]{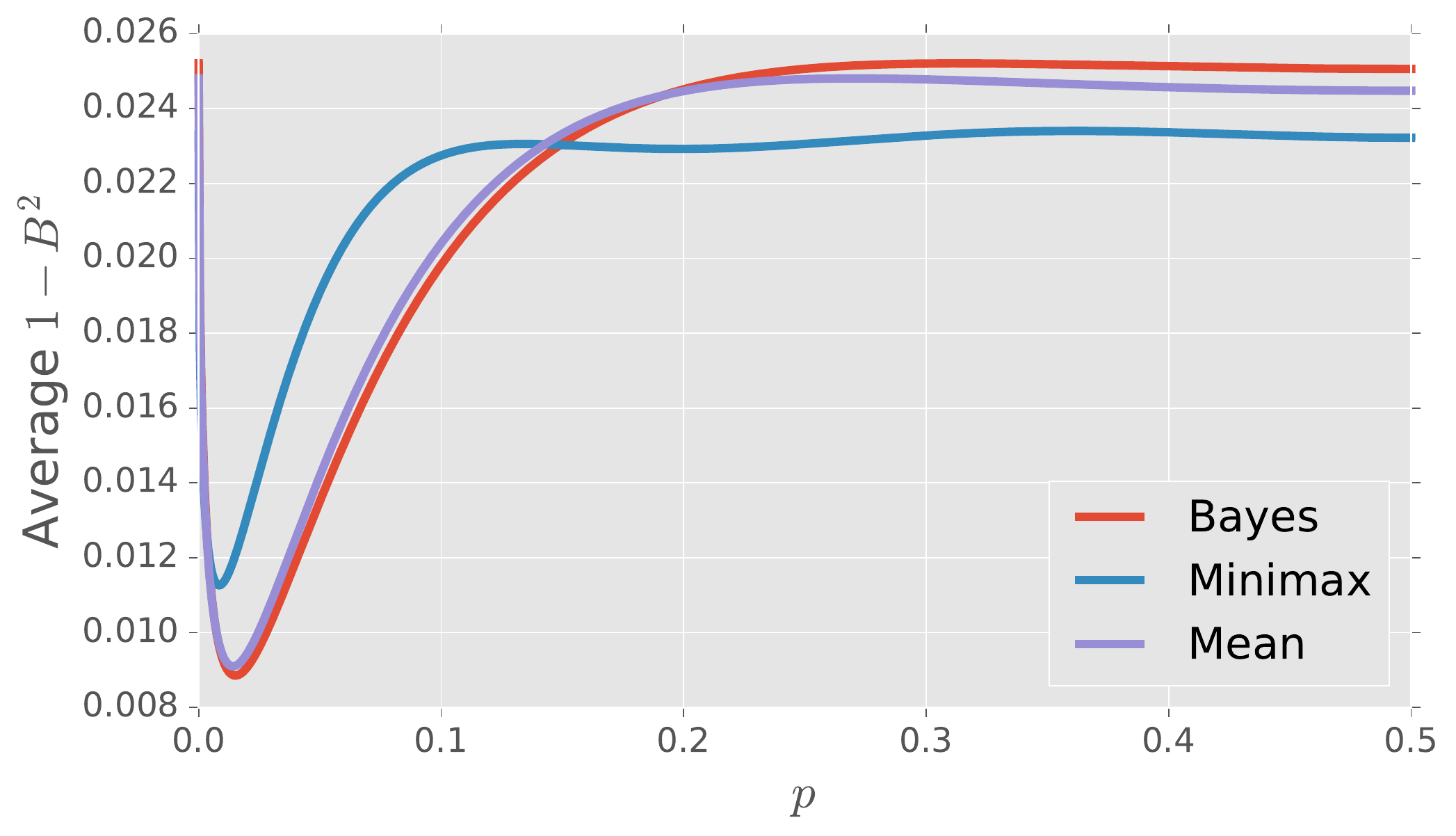}
\caption{\label{fig:minimax_N10}
In this figure, we compare three estimators:  (1) the minimax estimator, which is the Bayes estimator for the least favorable prior; (2) the Bayes estimator for the least favorable \emph{conjugate} prior, given by $\beta\approx0.44$; and (3) the posterior mean estimator for the least favorable [with respect to the posterior mean estimator] conjugate prior, given by $\beta\approx0.26$.  On the left, we show the estimators themselves.  On the right, we show their pointwise fidelity (1-risk) as a function of $p_0$.  $N=10$ in all cases.}
\end{figure}

To find the true minimax estimator, we used the algorithm of Kempthorne \cite{Kempthorne1987Numerical} to find the least favorable prior.  Roughly speaking, this algorithm iteratively adds support points to a discrete prior, each time maximizing the Bayes risk over the small finite search space of support points and weights.  It proceeds until the average risk and the maximum risk are within a pre-defined tolerance. We provide pseudo-code for our implementation of the algorithm in Algorithm \ref{kempthorne}. The results for the binomial parameter for $N=10$ are plotted in Figure \ref{fig:minimax_N10}.

\begin{algorithm}
\caption{\label{kempthorne}
Kempthorne's numerical algorithm for finding the least favorable prior \cite{Kempthorne1987Numerical}.}
\begin{algorithmic}
\Require Number of measurements $N>0$.
\Require Support points of initial guess $\vec p$.
\Require Probability weights of the support points $\vec w$.
\Require Tolerance $\mathtt{tol}>0$ (default $10^{-3}$).
\Require Mixing parameter $\alpha$ (default 0.01).
\Ensure   Least favorable prior $\{\vec{p},\vec{w}\}$.
\Ensure   Lower bound on the minimax risk $\mathtt{avg\_risk}$.
\Ensure   Upper bound on the minimax risk $\mathtt{max\_risk}$.
\Function{Kempthorne}{$N$, $\{\vec{p},\vec{w}\}$,$\mathtt{tol}$}
  \State $\mathtt{diff} \gets \mathtt{tol}+1$ \# make sure we enter the loop
  \While{$\mathtt{diff} > \mathtt{tol}$}
    \State $\{\vec{p},\vec{w}\} \gets$ prior with same number of support points which maximizes the Bayes risk
    \State $\mathtt{avg\_risk} \gets$ the maximum value of the Bayes risk for the prior found above
    \State $\mathtt{max\_risk} \gets$ global maximum of risk using the Bayes estimator of $\{\vec{p},\vec{w}\}$
    \State $\mathtt{diff}\gets |\mathtt{avg\_risk}-\mathtt{max\_risk}|/\mathtt{avg\_risk}$
    \If{$\mathtt{diff} > \mathtt{tol}$}
        \State Add a new support where the maximum risk is attained 
        \State $w_{\mathtt{length}(\vec{p})}\gets\alpha$ 
        \State for each $i\leq\mathtt{length}(\vec{p})-1$, $w_i\gets w_i-\alpha/(\mathtt{length}(\vec{p})-1)$   
    \EndIf 
  \EndWhile
  \State \Return $\{\vec{p},\vec{w}\},\mathtt{avg\_risk},\mathtt{max\_risk}$
\EndFunction
\end{algorithmic}
\end{algorithm}

\section{Conclusion\label{sec:6}}

Optimality of estimators is an always-relevant topic in statistics.  Even when optimal estimators are intractable or impractical, they provide a useful benchmark, and make it possible to show rigorously that some more tractable estimator is ``good enough''.  Our main technical contributions in this paper are simple, constructive formulas to compute Bayes estimators for ``fidelity''-type loss functions based on the Bhattacharyya coefficient.  This result may be directly useful for Bayesian machine learning, quantum state estimation, and other inference problems where it (1) provides a more accurate estimator, and (2) establishes a provable bound on performance.

We find the examples provided in the second half of our paper interesting because they suggest certain qualitative properties of multinomial estimation.  Two of the most widely used estimators for this problem are the MLE and the posterior mean (which is Bayes-optimal for Bregman divergences).  Our results indicate that the Bayes estimator for Bhattacharya loss interpolates between the MLE and the posterior mean; it hedges away from $p=0$ like the mean, but not as much.  Most interestingly, our analysis shows that all the estimators we considered have nearly identical average risk---and therefore suggests that worrying about finding exact Bayes estimators may be unnecessary.

Superficially, this may appear to undercut our result.  Who cares about deriving the Bayes estimator if it isn't significantly better?  But this is exactly the point:  our result \emph{proves} that the posterior mean is ``good enough'', at least in this particular case.  Furthermore, it shows that although the Bayes and mean \emph{estimators} are visibly different, their performance is not.  This suggests that---at least for Bhattacharya loss functions---a rather wide range of estimators achieve near-optimal performance.  This conclusion is reinforced by the behavior of the minimax (frequentist-optimal) estimators that we construct using our result, where we found that small changes in the precise definition of ``minimax'' produced fairly large changes in the ``optimal'' estimator.  And, while we did not fully solve the quantum problem (by finding Bayes estimators for quantum fidelity), we hope that our partial solution (for commuting states) provides a stepping stone to a full solution in the future.

\begin{acknowledgments}
CF acknowledges funding from the IARPA MQCO program, the ARC via EQuS project number CE11001013, and by the US Army Research Office grant numbers W911NF-14-1-0098 and W911NF-14-1-0103.  Sandia National Laboratories is a multi-mission laboratory managed and operated by Sandia Corporation, a wholly owned subsidiary of Lockheed Martin Corporation, for the U.S. Department of Energy's National Nuclear Security Administration 
under contract DE-AC04-94AL85000.
\end{acknowledgments}

\end{document}